\documentclass[12pt%,draft
]{article}

\usepackage[latin1]{inputenc}
\usepackage{amsfonts}
\usepackage{amsmath}
\usepackage{amssymb}
\usepackage{amscd}
\usepackage{amsthm}
\usepackage{indentfirst}
\usepackage[hmargin=3.6cm,vmargin=3.7cm]{geometry}

\usepackage{epsfig}

\newtheorem{thm}{Theorem}[section]
\newtheorem{lemma}[thm]{Lemma}
\newtheorem{prop}[thm]{Proposition}
\newtheorem{coroll}[thm]{Corollary}

\theoremstyle{definition}

\newtheorem{defin}[thm]{Definition}
\newtheorem{rem}[thm]{Remark}

\newtheorem{question}[thm]{Question}
\newtheorem*{acknow}{Acknowledgements}
\newtheorem*{prf}{Proof}

\newcommand{\R}{{\mathbb{R}}}

\newcommand{\T}{{\mathbb{T}}}

\newcommand{\N}{{\mathbb{N}}}
\newcommand{\C}{{\mathbb{C}}}

\newcommand{\cA}{{\mathcal{A}}}

\newcommand{\cH}{{\mathcal{H}}}

\newcommand{\cK}{{\mathcal{K}}}
\newcommand{\cL}{{\mathcal{L}}}

\newcommand{\cO}{{\mathcal{O}}}

\newcommand{\fc}{{:\ }}
\newcommand{\ve}{\varepsilon}

\newcommand{\ol}{\overline}

\newcommand{\tb}{\textbf}

\newcommand{\1}{\ensuremath 1\hspace{-0.25em}{\text{l}}}

\DeclareMathOperator{\id}{id}
\DeclareMathOperator{\area}{area}

\DeclareMathOperator{\pr}{pr}

\DeclareMathOperator{\Ham}{Ham}

\title{A comparison of symplectic homogenization and Calabi quasi-states}
\author{Alexandra Monzner\footnote{Fakult\"at f\"ur Mathematik, TU Dortmund; \texttt{alexandra.monzner@tu-dortmund.de}}\ { and} Frol Zapolsky\footnote{Max-Planck-Institut f\"ur Mathematik in den Naturwissenschaften; \texttt{zapolsky@mis.mpg.de}}}
\date{}

\begin{document}

\maketitle

%\begin{center}\textsc{Preliminary version}\end{center}

\begin{abstract}We compare two functionals defined on the space of continuous functions with compact support in an open neighborhood of the zero section of the cotangent bundle of a torus. One comes from Viterbo's symplectic homogenization, the other from the Calabi quasi-states due to Entov and Polterovich. In dimension $2$ we are able to say when these two functionals are equal. A partial result in higher dimensions is presented. We also give a link to asymptotic Hofer geometry on $T^*S^1$. Proofs are based on the theory of quasi-integrals and topological measures on locally compact spaces.
\end{abstract}

\renewcommand{\labelenumi}{(\roman{enumi})}

\section{Introduction and results}

A symplectic quasi-integral on a symplectic manifold is a positive functional on the space of continuous functions which is linear on Poisson commutative subspaces, and also satisfies a certain Lipschitz condition. Symplectic quasi-integrals (under the name of symplectic quasi-states) have been constructed on a variety of closed symplectic manifolds, for example, in \cite{EP_qs_sympl}, \cite{Ostrover_Calabi_qms_nonmonotone_mfds}, \cite{Usher_deformed_Floer}; applications can be found in \cite{EP_qs_sympl}, \cite{EPZ_qm_Poisson_br}. In the present paper we are interested in two particular examples of symplectic quasi-integrals. One comes from Viterbo's symplectic homogenization on $T^*\T^n$ \cite{Viterbo_homogenization}; the other is the Calabi quasi-state on $\C P^n$ due to Entov and Polterovich \cite{EP_qs_sympl}. Our goal is to compare the two. As a consequence of our computation of the homogenization operator in case $n=1$ we also obtain an explicit formula for the asymptotic Hofer norm of an autonomous Hamiltonian flow on $T^*S^1$ in terms of its homogenization.

%\subsection{Results}\label{section_results}

For a topological space $X$ we denote by $C(X)$ the space of real-valued continuous functions on $X$ while $C_c(X) \subset C(X)$ is the subspace of functions with compact support. We use the uniform norm $\|f\|=\sup_{x\in X}|f(x)|$ for $f\in C_c(X)$.
\begin{defin}\label{def_qi}
Let $X$ be a locally compact Hausdorff space. A (not necessarily linear) functional $\eta \fc C_c(X) \to \R$ is called a quasi-integral if
\begin{enumerate}
\item \tb{(Monotonicity)} $\eta(f)\leq\eta(g)$ for $f,g\in C_c(X)$ with $f \leq g$;
\item \tb{(Lipschitz continuity)} for every compact subset $K \subset X$ there is a number $N_K \geq 0$ such that $|\eta(f)-\eta(g)| \leq N_K\|f-g\|$ for all $f,g$ with support contained in $K$;
\item \tb{(Quasi-linearity)} $\eta$ is linear on every subspace of $C_c(X)$ of the form $\{\phi\circ f\,|\,\phi\in C(\R),\,\phi(0)=0\}$, where $f \in C_c(X)$.
\end{enumerate}
If $X$ is compact and $\eta$ satisfies $\eta(1)=1$, it is called a quasi-state; it is moreover called simple if $\eta(f^2)=\big(\eta(f)\big)^2$ for all $f\in C(X)$. If $X$ is a symplectic manifold and $\eta$ is linear on Poisson commutative subspaces of $C_c^\infty(X)$, it is called symplectic.
\end{defin}

\begin{rem} In case $X$ is compact, the notion of a quasi-integral was introduced and first studied by Aarnes, see \cite{Aarnes_quasi-states}. It was generalized to various other settings, see Remark \ref{rem_other_generalizations_qis} for more information. However, as far as we know, the notion of a quasi-integral as a positive quasi-linear Lipschitz functional on the space of continuous functions with compact support on a locally compact space is new.
\end{rem}

Quasi-integrals are a generalization of integration against a Radon measure; the latter are exactly the linear quasi-integrals.

There is a representation theory (extending the Riesz representation theorem) for quasi-integrals in terms of certain set functions, called topological measures, which we describe in detail below. For now let us mention that the value of the topological measure representing a quasi-integral $\eta$ on a compact subset $K \subset X$ equals, intuitively speaking, the value of $\eta$ on the (discontinuous) indicator function of $K$.

\subsection{Dimension two}\label{section_dim_two}

First we are going to introduce the relevant functionals in dimension two. They can be uniquely characterized by simple properties. The first functional is the unique simple quasi-state $\zeta$ on $S^2$ which is invariant under Hamiltonian diffeomorphisms (relative to the standard symplectic structure), the so-called Calabi quasi-state, see \cite{EP_qs_sympl}. Identify $T^*S^1 = S^1 \times \R$; then we also have the following lemma, whose proof is given in subsection \ref{section_computation}:
\begin{lemma}\label{lemma_unique_homogenization}There is a unique (non-linear) operator $\cH \fc C_c(T^*S^1) \to C_c(\R)$ such that for all $f,g \in C_c(T^*S^1)$:

(i) \tb{(Monotonicity)} $\cH(f)\leq\cH(g)$ if $f \leq g$;

(ii) \tb{(Lipschitz continuity)} $\|\cH(f)-\cH(g)\| \leq \|f-g\|$;

(iii) \tb{(Strong quasi-linearity)} the restriction of $\cH$ to any Poisson commu\-ta\-ti\-ve subspace of $C^\infty_c(T^*S^1)$ is linear;

(iv) \tb{(Invariance)} if $\phi$ is a Hamiltonian diffeomorphism of $T^*S^1$ generated by a time-dependent Hamiltonian with compact support, then $\cH(f\circ\phi) = \cH(f)$;

(v) \tb{(Lagrangian)} if there is a constant $c$ and $p\in \R$ so that $f = c$ on $S^1\times\{p\}$, then $\cH(f)(p)=c$.
\end{lemma}

\noindent Define $\eta_0 \fc C_c(T^*S^1) \to \R$ by $\eta_0(f)=\cH(f)(0)$. For $r \in (0,\frac 1 2]$ consider a symplectic embedding $j \fc U_r \to S^2$, where $U_r = S^1\times (-r,r) \subset T^*S^1$, such that $j(S^1 \times \{0\})$ is the equator. The symplectic forms are standard and they are normalized so that the area of $S^2$ is $1$ while the area of $U_r$ is $2r$. There is the induced extension-by-zero map $j_! \fc C_c(U_r) \to C(S^2)$ and the pull-back functional $\zeta_r:=j^*\zeta = \zeta \circ j_!$ on $C_c(U_r)$. The following is a comparison between the functionals $\zeta_r$ and $\eta_0$:
\begin{thm}\label{thm_two_dim}The restriction of $\eta_0$ to $C_c(U_r)$ coincides with $\zeta_r$ if and only if $r \in (0,\frac 1 4]$.
\end{thm}

The functional $\zeta$ and the operator $\cH$ are particular cases of more general constructions, which we describe in the following subsection.

\subsection{The functionals}\label{section_functionals}

\subsubsection{The Calabi quasi-state}\label{section_Calabi_qs}

The Calabi quasi-state on $\C P^n$ is the stabilization of a certain Hamiltonian Floer-homological spectral invariant. In more detail, consider $\C P^n$ with its standard symplectic structure $\omega$, normalized so that $\int_{\C P^n}\omega^n = 1$. We refer the reader to \cite{Hofer_Salamon_Floer_hlgy_Novikov_rings}, \cite{McDuff_Salamon_J_hol_curves}, \cite{Oh_spectral_invariants} for preliminaries on Floer and quantum homology with coefficients in a Novikov ring, as well as on Oh's spectral invariants of Hamiltonian diffeomorphisms. The quantum homology $QH$ of $\C P^n$ with coefficients in a suitable Novikov ring is an associative algebra with respect to the quantum product, which is indeed a field. Let $e \in QH$ denote the unit element of this field, given by the fundamental class $[\C P^n]$. We let $c \fc QH-\{0\} \times \widetilde{\Ham}(\C P^n) \to \R$ be the spectral invariant defined on the universal cover $\widetilde \Ham (\C P^n)$ of the group of Hamiltonian diffeomorphisms of $\C P^n$. The Calabi quasi-state $\zeta$ on $\C P^n$ is defined as follows. First, for $f \in C^\infty(\C P^n)$, let $\phi_f$ denote the element of $\widetilde \Ham$ generated by $f$. Then
$$\zeta(f)=\lim_{k \to \infty}\frac {c(e,\phi_f^k)}{k} + \int_{\C P^n}f\,\omega^n\,.$$
This functional is Lipschitz with respect to the $C^0$ norm and so admits a unique extension to $C(\C P^n)$, which is a symplectic quasi-state, called the Calabi quasi-state \cite{EP_qs_sympl}.

Let $\tau$ denote the topological measure representing the Calabi quasi-state. It is proved in \cite{EP_qs_sympl} that the Calabi quasi-state is invariant under Hamiltonian diffeomorphisms. The same is then true for $\tau$. We also need the following results about $\tau$, proved \emph{ibid.}

\emph{The case $n=1$.} Here $\C P^1 = S^2$, with the area form normalized to have area $1$. A topological measure is determined by its values on compact connected subsurfaces with boundary (Lemma \ref{lemma_tms_uniquely_detd_by_sbmfds}). Any such subsurface of $S^2$ has the form $W=S^2- \bigcup_i D_i$, where the $D_i$ are finitely many open disks with disjoint closures. We have (see for example \cite{Aarnes_qm_construct} for a proof that this indeed defines a topological measure)
$$\tau(W)=\left\{\begin{array}{ll}0\,,&\text{if }\area(D_i) > \frac 1 2 \text{ for some }i\\ 1\,,&\text{otherwise}\end{array}\right.\;.$$
It also follows that if $D$ is an open disk, then $\tau(D) = 0$ if $\area(D) \leq \frac 1 2$ and $1$ otherwise. Finally, if $L\subset S^2$ is an equator, that is, a simple closed curve such that its complement is two open disks of area $\frac 1 2$, then $\tau(L)=1$.

\emph{The case $n\geq 2$.} Let $\T^n_{\text{Clif}} = \{[z_0:\dots:z_n]\,|\,|z_0|=\dots=|z_n|\} \subset \C P^n$ be the Clifford torus. Then $\tau(\T^n_{\text{Clif}}) = 1$. Since $\tau$ is invariant under Hamiltonian isotopies, if $L$ is Hamiltonian isotopic to $\T^n_{\text{Clif}}$, then $\tau(L) = 1$.

\subsubsection{The symplectic homogenization}

The general reference for what appears in this subsection is the paper of Viterbo \cite{Viterbo_homogenization}. The homogenization $\cH(f) \fc \R^n \to \R$ of $f \in C_c(T^*\T^n)$ is the limit with respect to Viterbo's metric of the sequence $f_k(q,p)=f(kq,p)$, $k \in \N$. More precisely, a Hamiltonian isotopy $\phi_t$ generated by a compactly supported time-dependent Hamiltonian yields the Hamiltonian isotopy $\id \times \phi_t$ of $T^*\T^n \times \ol{T^*\T^n}$, where the overline indicates that the symplectic form has the opposite sign. This latter manifold has a symplectic covering by $T^*(\Delta_{T^*\T^n})$, $\Delta$ being the diagonal in the product, therefore there is a unique Hamiltonian lift $\widetilde \phi_t$ such that $\widetilde \phi_0 = \id$, and we consider the Lagrangian submanifold $L=\widetilde \phi_1(O_{T^*\T^n})$, $O$ denoting the zero section. $L$  has an essentially unique generating function quadratic at infinity $S \fc T^*\T^n \times E \to \R$ ($E$ is a parameter vector space), normalized to be equal to zero at the infinity of $T^*\T^n$ \cite{Viterbo_gfs}. For $p \in \R^n$ denote $S_p(q,\xi) = S(q,p,\xi)$ and let $h(p) = c(\mu_{\T^n},S_p)$. Here $c$ is the Lagrangian spectral invariant due to Viterbo \cite{Viterbo_gfs} and $\mu_{\T^n} \in H^n(\T^n)$ is the orientation class. Now given $f \in C^\infty_c(T^*\T^n)$ define, for $k \in \N$, $f_k(q,p)=f(kq,p)$, consider the time-$1$ flow $\phi_k$ of $f_k$, and let $h_k \fc \R^n \to \R$ denote the function $h$ associated to $\phi_k$ by the above construction. The homogenization $\cH(f)$ of $f$ is then the limit of a subsequence $\{h_{k_m}\}_m$ in the $C^0$ norm. Moreover, the time-$1$ maps $\phi_k$ form a Cauchy sequence with respect to Viterbo's metric and its limit in the completion with respect to this metric is in a certain precise sense generated by the $q$-independent Hamiltonian $\cH(f)$ \cite{Viterbo_homogenization}.

The resulting operator $\cH \fc C^\infty_c(T^*\T^n) \to C_c(\R^n)$ is Lipschitz with respect to the $C^0$ norm and so admits a unique extension to $C_c(T^*\T^n)$, which is the symplectic homogenization operator. This operator has the properties analogous to those described in Lemma \ref{lemma_unique_homogenization}, with $T^*S^1$ replaced by $T^*\T^n$. We again denote $\eta_0(f)=\cH(f)(0)$. For future use we formulate
\begin{lemma}\label{thm_Viterbo_fcnls_are_qis}Let $\mu$ be a Radon measure\footnote{A Radon measure for us is a locally finite regular Borel measure.} on $\R^n$. Then
$$\eta_\mu :=\int_{\R^n}\cH(\cdot)\,d\mu \fc C_c(T^*\T^n) \to \R$$
is a symplectic quasi-integral; it is Lipschitz continuous with constant $\mu(K)$ for functions with compact support in $\T^n \times K$, where $K \subset \R^n$ is compact.
\end{lemma}
\noindent It follows that $\eta_0 = \eta_{\delta_0}$ is a symplectic quasi-integral.
\begin{rem}In \cite{Viterbo_homogenization}, Viterbo formulated a version of this lemma, and proved most of what is stated in it, although the proper definition of a quasi-integral on $T^*\T^n$ was lacking.
\end{rem}

\subsection{Higher dimensions}\label{section_intro_higher_dims}

Our goal is a comparison of Calabi quasi-states and symplectic homogenization. Therefore next we describe a partial negative result in dimensions $n \geq 2$. Let us first give some context and motivation. Both $\C P^n$ and $T^*\T^n$ admit Hamiltonian torus actions, with moment maps $\Phi \fc \C P^n \to \R^n$ and $\Psi \fc T^*\T^n \to \R^n$ defined by
$$\Phi([z_0:\dots:z_n])=\left(\frac {|z_1|^2}{\sum_{j=0}^n |z_j|^2}-\frac 1 {n+1},\dots,\frac {|z_n|^2}{\sum_{j=0}^n |z_j|^2}-\frac 1 {n+1}\right)$$
and
$$\Psi(q,p)=p\,,$$
where we view $T^*\T^n = \T^n(q) \times \R^n(p)$. The functionals $\zeta$ and $\eta_0$ satisfy
$$\zeta(\Phi^*\ol f)=\ol f(0),\qquad \eta_0(\Psi^* \ol g) = \ol g(0)$$
for $\ol f \in C(\R^n)$ and $\ol g \in C_c(\R^n)$. For $\eta_0$ this follows from the `Lagrangian' property of $\cH$; for $\zeta$ this is proved in \cite{EP_qs_sympl}. Another way of saying this is as follows. Symplectic quasi-integrals push forward by proper maps to quasi-integrals. If the map in question is to $\R^n$ and its coordinate functions Poisson commute, the resulting quasi-integral is globally linear and so corresponds to a measure. The above can be restated by saying that the push-forwards $\Phi_*\zeta$ and $\Psi_*\eta_0$ are both delta-measures at $0 \in \R^n$, since the coordinates of a moment map of a Hamiltonian torus action Poisson commute.

There exists a symplectic embedding with dense image $j \fc U \to \C P^n$, where $U =\T^n \times V \subset T^*\T^n$, $V = \{p \in \R^n\,|\,p_j > -\frac 1 {n+1},\, \sum_j p_j < \frac 1 {n+1}\}$, which commutes with moment maps, that is $\Phi \circ j = \Psi$. Consider again the induced map $j_! \fc C_c(U) \to C(\C P^n)$ and the pull-back functional $j^*\zeta = \zeta \circ j_!$. It follows that if $F = \Psi^*\ol f$ for $\ol f \in C_c(V)$, then $j^*\zeta(F) = \eta_0(F)$, that is, the two functionals agree on functions pulled back via the moment map. It is then natural to pose the following question:
\begin{question}[L.\ Polterovich] Do $j^*\zeta$ and $\eta_0$ agree on all of $C_c(U)$? If not, does there exist an open neighborhood $V' \subset \R^n$ of $0$ such that they do on $C_c(U')$ where $U' = \T^n \times V'$?
\end{question}

Theorem \ref{thm_two_dim} answers this question for $n=1$. A partial negative result for $n \geq 2$ is:
\begin{prop}\label{prop_counterex_higher_dim}For any $\ve, \delta \in (0,\frac 1 {n(n+1)})$ consider $U_{\ve,\delta} = \T^n \times V_{\ve,\delta}$, where
$$\textstyle V_{\ve,\delta}= (-\frac 1{n+1},\ve) \times (-\delta,\delta)^{n-1} \subset \R^n\,.$$
Then the restrictions of $\eta_0$ and $j^*\zeta$ to $C_c(U_{\ve,\delta})$ do not coincide.
\end{prop}

\begin{rem}Note that Theorem \ref{thm_two_dim} is a comparison of the functional $\zeta_r$ coming from the Calabi quasi-state and the functional $\eta_0$ in dimension $2$. For measures $\mu$ (not necessarily delta-measures) with support in $(-r,r)$, other than the delta-measure at $0$, the functionals $\zeta_r$ and $\eta_\mu$ \emph{do not} coincide, as can be seen already by evaluating the two on functions pulled back by the moment map $\Psi$. A similar remark applies in higher dimensions. Therefore throughout we only speak about $\eta_0$.
\end{rem}

\begin{rem}\label{rem_uniqueness}We draw the reader's attention to the fact that the above results have to do with the general question of uniqueness of symplectic quasi-states and quasi-integrals. Linear quasi-integrals are in one-to-one correspondence with Radon measures. In general, a positive linear combination of a non-linear symplectic integral and of a linear one yields a non-linear one, so the interesting question is whether symplectic integrals are unique up to the addition of a measure. As Theorem \ref{thm_two_dim} and Proposition \ref{prop_counterex_higher_dim} show, there is no uniqueness of symplectic integrals on a neighborhood of the zero section in $T^*\T^n$, even if we impose additional conditions like Hamiltonian invariance (compare with Lemma \ref{lemma_unique_homogenization}) and the values of the quasi-integral on functions pulled back from the moment map $\Psi$.
\end{rem}

\subsection{A link to asymptotic Hofer geometry}

Given a symplectic manifold $(M,\omega)$ and a Hamiltonian diffeomorphism $\phi$ of $M$ generated by a time-dependent Hamiltonian with compact support, the Hofer norm\footnote{See \cite{Polterovich_geom_grp_sympl_diffeo} for preliminaries on Hofer geometry.} of $\phi$ is the number
$$\|\phi\|_{\text{Hofer}} = \inf_f \int_0^1\Big(\max_M f_t - \min_M f_t\Big)\,dt\,,$$
where the infimum is taken over all compactly supported Hamiltonians $f\fc[0,1]\times M \to \R$ generating $\phi$ and $f_t(\cdot) \equiv f(t,\cdot)$. For an autonomous Hamiltonian $f$ on $M$ denote by $\phi_f^t$ its flow and define the asymptotic Hofer norm of (the flow of) $f$ by
$$\mu_\infty(f)=\lim_{t\to+\infty}\frac{\|\phi_f^t\|_{\text{Hofer}}} t\,.$$
We refer the reader to \cite{Polterovich_Siburg_asym_geom_area_prsv_maps} for a discussion on $\mu$. \emph{Ibid.}, the authors prove that in case $M$ is an open surface of infinite area, it is true that
$$\mu_\infty(f)=c_+(f)-c_-(f)\,,$$
where
$$c_+(f)=\sup_{L \in \cL} \min_L f \qquad \text {and} \qquad c_-(f) = \inf_{L \in \cL}\max_L f\,,$$
$\cL$ being the set of all embedded non-contractible circles in $M$. A stronger result proved there is that $t\mapsto\|\phi^t_f\|_{\text{Hofer}}$ is either bounded or asymptotically linear, depending on whether $c_+$ equals $c_-$. We show the following
\begin{lemma}\label{lemma_computation_of_c_pm}If $M=T^*S^1$ with its canonical symplectic structure, then for $f \in C^\infty_c(T^*S^1)$ we have $c_+(f) = \max \cH(f)$ and $c_-(f) = \min \cH(f)$.
\end{lemma}
\noindent The proof is given in subsection \ref{section_eta_p_prf_lemma}. This implies
\begin{coroll}The asymptotic Hofer norm of $f \in C^\infty_c(T^*S^1)$ satisfies
$$\mu_\infty(f)=\max \cH(f)-\min\cH(f)\,.\qed$$
\end{coroll}

In \cite{Viterbo_homogenization} it is proven that the quantity $\max \cH(f)-\min\cH(f)$ also equals the so-called asymptotic Viterbo distance $\gamma_\infty(f)$, see \cite{Sorrentino_Viterbo_act_min_pties_dist_Ham} for definitions. It is true in general that the (asymptotic) Viterbo distance is bounded from above by the (asymptotic) Hofer norm (\emph{ibid.}). The above discussion shows that in the autonomous case we have the following corollary, where $\gamma$ is the usual Viterbo distance:
\begin{coroll}For $f \in C^\infty_c(T^*S^1)$ we have $\gamma_\infty(f)=\mu_\infty(f)$, moreover, $t\mapsto\gamma(\phi^t_f)$ is either bounded or asymptotically linear. \qed
\end{coroll}

\begin{rem}This should be contrasted with Theorem 2 of the same paper which gives a construction of a fiberwise convex autonomous Hamiltonian on the closed disk cotangent bundle $B^*\T^n$ of a torus for which $\gamma_\infty$ is strictly less than the asymptotic Hofer norm. The latter however is defined using only Hamiltonians on $B^*\T^n$ which vanish on the boundary and which admit a smooth extension to $T^*\T^n$ depending only on time and on $\|p\|$ outside $B^*\T^n$ (this particular flavor of asymptotic Hofer geometry was introduced and first studied in \cite{Siburg_act_min_meas_geom_Ham}). The fact that $B^*\T^n$ has finite volume allows to use the Calabi invariant of the Hamiltonian as a lower bound for its asymptotic Hofer norm, which is impossible on $T^*S^1$.
\end{rem}

\subsection{Topological measures}

In order to compare $\zeta$ and $\eta_0$, we make use of a representation theorem for quasi-integrals in terms of topological measures. %The reader is referred to section \ref{section_qis_tms_loc_cpt} for a detailed treatment.

\begin{defin}\label{def_tm}Let $X$ be a locally compact Hausdorff space. Let $\cK(X)$ be the family of compact subsets of $X$, $\cO(X)$ the family of open subsets of $X$ with compact closure, and $\cA(X) = \cK(X) \cup \cO(X)$. A function $\tau \fc \cA(X) \to [0,\infty)$ is a topological measure if
\begin{enumerate}
\item \tb{(Additivity)} if $A,A'\in\cA(X)$ are disjoint and $A\cup A'\in \cA(X)$, then $\tau(A\cup A') = \tau(A)+\tau(A')$;
\item \tb{(Monotonicity)} $\tau(A)\leq \tau(A')$ for $A,A' \in \cA(X)$ with $A\subset A'$;
\item \tb{(Regularity)} $\tau(K)=\inf\{\tau(O)\,|\,O\in\cO(X),\,O\supset K\}$ for any $K \in \cK(X)$ (outer) and $\tau(O)=\sup\{\tau(K)\,|\,K\in\cK(X),\,O\supset K\}$ for any $O \in \cO(X)$ (inner).
\end{enumerate}
\end{defin}

\begin{rem} A topological measure in this sense on a compact space is the same as a usual Aarnes topological measure \cite{Aarnes_quasi-states}.
\end{rem}

\begin{thm}There is a natural bijection between the sets of quasi-integrals and of topological measures on a locally compact space.
\end{thm}

The description of the bijection, as well as a more precise formulation of the theorem and its proof are the subject of section \ref{section_qis_tms_loc_cpt}.

\begin{acknow}We thank Leonid Polterovich for suggesting the topic of this paper and for his interest in it, as well as for pointing out the link to asymptotic Hofer geometry, and Karl Friedrich Siburg for useful discussions and suggestions. The second author would like to thank Judy Kupferman and Marco Mazzucchelli for listening to a preliminary version of the results and for helpful suggestions. This work started during the stay of the first author at the University of Chicago, which was supported by the Martin-Schmeißer-Foundation. The first author is partially supported by the German National Academic Foundation. This work is partially supported by the NSF-grant DMS 1006610.
\end{acknow}

\section{Quasi-integrals and topological measures on locally compact spaces}\label{section_qis_tms_loc_cpt}

In this section $X$ is a locally compact Hausdorff space. Any open subset of $X$ is a locally compact space on its own right. Note also that $X$ is completely regular; we will implicitly use this fact and its consequences. We use the theory for the compact case, developed in \cite{Aarnes_quasi-states}, without explicitly mentioning it.

Recall the definitions of a quasi-integral and of a topological measure on $X$, Definitions \ref{def_qi}, \ref{def_tm}.

Given a quasi-integral $\zeta$, define a set function $\tau_\zeta \fc \cA(X) \to [0,\infty)$ by
$$\tau_\zeta(K) = \inf \{\zeta(f)\,|\,f\in C_c(X),\,f\geq\1_K\};\,\tau_\zeta(O) = \sup \{\zeta(f)\,|\,f\in C_c(X),\,f\leq\1_O\}$$
for $K \in \cK(X)$ and $O\in\cO(X)$. Here and in the sequel, $\1$ stands for the characteristic function of a set.

The main result of this section is
\begin{thm}[Representation Theorem]\label{thm_representation}The map $\zeta \mapsto \tau_\zeta$ is a bijection from the space of quasi-integrals to the space of topological measures on $X$.
\end{thm}
Most of this section is devoted to the proof of this theorem. We would like to point out that in the original work \cite{Aarnes_quasi-states} Aarnes used delicate analysis in order to prove his representation theorem. Instead of adapting his arguments to the locally compact case, we rely on results valid in the compact case, the rest of the proof being relatively elementary.

Subsection \ref{section_one_pt_cpct} contains the main technical step which allows a reduction to the compact case. Subsection \ref{section_prf_repr_thm} is devoted to the proof of the representation theorem. Subsection \ref{section_prf_Viterbo_fcnls_are_qis} contains the proof of Lemma \ref{thm_Viterbo_fcnls_are_qis}.

\subsection{One-point compactifications}\label{section_one_pt_cpct}

Fix $O \in \cO(X)$ and let $\widehat O = O \cup \infty$ be its one-point compactification. Recall that
$$\cO(\widehat O)=\{U\subset O\text{ open}\} \cup \{(O-K)\cup \infty\,|\,K\in\cK(O)\}$$
and
$$\cK(\widehat O) = \cK(O)\cup \{(O-U)\cup \infty\,|\, U \subset O \text{ open}\}\,.$$

Fix a topological measure $\tau$ on $X$ and define $\widehat \tau_O \fc \cA(\widehat O) \to [0,\infty)$ by
$$\widehat \tau_O(U)=\tau(U)\,,\quad \widehat\tau_O(K) = \tau(K)$$
and
$$\widehat \tau_O((O-K)\cup \infty)= \tau(O-K)\,,\quad \widehat \tau_O((O-U)\cup\infty)=\tau(O)-\tau(U)$$
for $U\subset O$ open and $K\in\cK(O)$.

\begin{lemma}\label{lemma_tm_on_compactification}$\widehat \tau_O$ is a topological measure on $\widehat O$.
\end{lemma}

The proof is a routine verification; we supply it for the sake of completeness.

\begin{prf}Note first of all that $\widehat \tau \equiv \widehat \tau_O$ is well-defined. We need to show that
\begin{enumerate}
\item $\widehat\tau(\widehat O - K) + \widehat \tau(K) = \widehat\tau(\widehat O)$ for $K \in\cK(\widehat O)$;
\item $\widehat\tau(K \cup K')=\widehat\tau(K) + \widehat\tau(K')$ for disjoint $K,K' \in \cK(\widehat O)$;
\item $\widehat\tau(K)\leq\widehat\tau(K')$ for $K,K' \in \cK(\widehat O)$ with $K\subset K'$;
\item $\widehat\tau(K)=\inf\{\widehat\tau(U)\,|\,U\in\cO(\widehat O),\,U \supset K\}$ for $K\in\cK(\widehat O)$.
\end{enumerate}

We note the following: if $K,K'$ are compact subsets of $O$, then all of the above properties follow immediately from the definition of $\widehat \tau$ and the corresponding properties of $\tau$. The following then suffices to establish (i-iv). (i) Let $K=(O-U)\cup \infty$, where $U\subset O$ is open. Then
$$\widehat \tau(\widehat O - K)+\widehat\tau(K)=\tau(U)+(\tau(O)-\tau(U))=\tau(O)=\widehat\tau(\widehat O)\,.$$
(ii) Let $K\in\cK(O)$ and $K' = (O-U)\cup\infty\in\cK(\widehat O)$ be disjoint, where $U \subset O$ is open. Then
$$\widehat\tau(K\cup K')=\widehat\tau(O-(U-K)\cup\infty)=\tau(O)-\tau(U-K)=\underbrace{(\tau(O)-\tau(U))}_{=\widehat\tau(K')}+\underbrace{\tau(K)}_{=\widehat\tau(K)}\,.$$
(iii) Let $K\in\cK(O)$ and $K' = (O-U)\cup\infty\in\cK(\widehat O)$ be such that $K\subset K'$, where $U \subset O$ is open. Then $K\subset O-U$, which implies $U \subset O-K$, hence $\tau(U) \leq \tau(O-K)=\tau(O)-\tau(K)$, and we have
$$\widehat\tau(K')=\tau(O)-\tau(U) \geq \tau(K) = \widehat \tau(K)\,.$$
Now assume $K = (O-V)\cup \infty \in \cK(\widehat O)$, with $V \subset O$ open and $K \subset K'$, where $K'$ is as above. Then $V \supset U$ and so
$$\widehat\tau(K)=\tau(O)-\tau(V) \leq \tau(O)-\tau(U)=\widehat\tau(K')\,.$$
(iv) Let $K = (O-U)\cup\infty\in\cK(\widehat O)$, where $U \subset O$ is open. If an opet set $V\subset \widehat O$ contains $K$, it has to be of the form $V=(O-L)\cup \infty$, where $L\subset O$ is compact, and then it follows that $L \subset U$. Thus
$$\inf\{\widehat\tau(V)\,|\,V\supset K\text{ open}\} = \inf\{\tau(O)-\tau(L)\,|\,L\subset U\text{ compact}\}\,,$$
which equals
$$\tau(O)-\sup\{\tau(L)\,|\,L\subset U\text{ compact}\}=\tau(O)-\tau(U)=\widehat\tau(K)\,,$$
where the first equality follows from the inner regularity of $\tau$.
\qed
\end{prf}

Now we apply the one-point compactification process to quasi-integrals. Let $\zeta$ be a quasi-integral. The space $C_c(O)$ is dense in the $C^0$ norm in the space $C_0(\widehat O)=\{f\in C(\widehat O)\,|\,f(\infty)=0\}$. The restriction of $\zeta$ to $C_c(O)$ is Lipschitz and so defines a unique extension $\widehat\zeta_O$ to $C_0(\widehat O)$, which is also Lipschitz, with the same Lipschitz constant as $\zeta|_{C_c(O)}$. For a general $f\in C(\widehat O)$ put $\widehat \zeta_O(f)=\widehat \zeta_O(f-f(\infty))+\lambda_O\cdot f(\infty)$, where $\lambda_O = \sup \{\zeta(g)\,|\,g\in C_c(O),\,g \leq \1_O\}$. Although $\lambda_O = \tau_\zeta(O)$, we will ignore this fact for the moment, because we need to express everything in terms of $\zeta$.
\begin{lemma} $\widehat \zeta_O$ is a quasi-integral on $\widehat O$.
\end{lemma}
\begin{prf}Abbreviate $\widehat\zeta = \widehat \zeta_O$. It suffices to show (i) $\widehat\zeta(f) \geq 0$ for $f \in C(\widehat O)$, $f\geq 0$, and (ii) $\widehat\zeta$ is linear on every subspace of $C(\widehat O)$ of the form $\{\phi\circ f\,|\,\phi \in C(\R)\}$ with $f\in C(\widehat O)$.

Proof of (i). Let $f \in C(\widehat O)$, $f \geq 0$. Put $\widetilde f = f - f(\infty)$. We have $\widehat\zeta(f)=\widehat\zeta(\widetilde f) + \lambda_O f(\infty)$. Let $\ve > 0$. There is $g \in C_c(O)$ such that $\|\widetilde f-g\| < \ve$ and $0 \geq \min g = \min \widetilde f \geq -f(\infty)$. There exists $h \in C_c(O)$ such that $0 \geq h \geq \min g$ and $h = \min g$ on the support of $g$. We then have $g \geq h \geq \min g\cdot \1_O$ and so
$$\widehat\zeta(g) = \zeta(g) \geq \zeta(h) \geq \lambda_O \cdot\min g = \lambda_O\cdot \min \widetilde f \geq -\lambda_O\cdot f(\infty)\,,$$
by the definition of $\lambda_O$ and the linearity of $\zeta$ on $\R\cdot h \subset C_c(O)$. Therefore
$$\widehat\zeta(\widetilde f) \geq \widehat\zeta(g) - C\ve \geq \zeta(h) - C\ve \geq -\lambda_O\cdot f(\infty) - C\ve\,,$$
where the first inequality follows from the Lipschitz continuity of $\widehat\zeta|_{C_0(\widehat O)}$ with constant $C$. Thus we obtained, for any $\ve >0$:
$$\widehat\zeta(f)=\widehat\zeta(\widetilde f) + \lambda_O \cdot f(\infty) \geq -C\ve\,,$$
which proves (i).

For (ii) let $f \in C(\widehat O)$ and $\phi,\psi \in C(\R)$. Since
$$\widehat\zeta(\phi\circ f)=\widehat\zeta(\phi\circ f - \phi(f(\infty))) + \lambda_O\phi(f(\infty)) = \widehat\zeta((\phi - \phi(f(\infty)))\circ f) + \lambda_O\phi(f(\infty))$$
and similarly for $\psi\circ f$ and $\phi\circ f+\psi\circ f = (\phi+\psi)\circ f$, proving that $\widehat\zeta(\phi\circ f + \psi\circ f) = \widehat\zeta(\phi \circ f) + \widehat\zeta(\psi\circ f)$ is equivalent to proving that $\widehat\zeta[(\phi - \phi(f(\infty)))\circ f + (\psi - \psi(f(\infty)))\circ f]=\widehat\zeta((\phi- \phi(f(\infty)))\circ f)+\widehat\zeta((\psi - \psi(f(\infty)))\circ f)$, therefore we may assume that $\phi(f(\infty))=\psi(f(\infty))=0$. Now let $\widetilde f = f-f(\infty)$, and let $\widetilde \phi$, $\widetilde \psi$ be defined by $\widetilde \phi(t)=\phi(t+f(\infty))$ and similarly for $\widetilde \psi$. We then have $\widetilde \phi(0) = \widetilde \psi(0)=0$ and $\widetilde \phi \circ \widetilde f = \phi \circ f$ and same for $\psi$. Let now $f_k \in C_c(O)$ be a sequence whose limit is $\widetilde f$. The fact that $\zeta$ is quasi-linear implies that
$$\zeta(\widetilde \phi \circ f_k + \widetilde \psi \circ f_k)=\zeta(\widetilde \phi \circ f_k)+\zeta(\widetilde \psi \circ f_k)\,.$$
When $k \to \infty$, the left-hand side tends to $\widehat\zeta(\widetilde \phi \circ \widetilde f + \widetilde \psi \circ \widetilde f) = \widehat\zeta(\phi \circ f + \psi \circ f)$, while the right-hand side to $\widehat\zeta(\widetilde \phi \circ \widetilde f) + \widehat\zeta(\widetilde \psi \circ \widetilde f) = \widehat\zeta(\phi \circ f) + \widehat\zeta(\psi \circ f)$, thereby proving (ii). \qed

\end{prf}

\subsection{Proof of Theorem \ref{thm_representation}}\label{section_prf_repr_thm}

\subsubsection{From quasi-integrals to topological measures}

Recall that we defined a set function $\tau_\zeta$ using a quasi-integral $\zeta$.

\begin{prop}$\tau=\tau_\zeta$ is a topological measure.
\end{prop}

\begin{prf}
\tb{Monotonicity}: for pairs of compact subsets, as well as for pairs of open subsets follows from the definition. If $K \in \cK(X)$, $O \in \cO(X)$ and $O \subset K$, then for any function $f$ such that $f \geq \1_K$ and any function $g$ with $g \leq \1_O$ we have $f \geq g$ and so $\tau(K) = \inf \zeta(f) \geq \sup \zeta(g) = \tau(O)$, the $\inf$ and $\sup$ being taken over all such $f,g$. Assume now that $K \subset O$. Then there exist $f \in C_c(X)$ with values in $[0,1]$ such that $f|_K =1$ and $f|_{X-O}=0$. Thus $\tau(K) \leq \zeta(f) \leq \tau(O)$.

\tb{Regularity}: let $K\in\cK(X)$. For outer regularity we have to prove that $\tau(K) = \inf \{\tau(O)\,|\,O\in\cO(X),\,O\supset K\}$. Denote the infimum by $I$. Then from monotonicity it follows that $\tau(K) \leq I$ and we want to show that $\tau(K) \geq I$. Let $\ve > 0$, and fix a compact set $L$ containing $K$ in its interior. By the definition of the infimum and the fact that $X$ is completely regular, there is a function $f$ such that $f|_K =1$, $f=0$ outside the interior of $L$, and $\tau(K) \geq \zeta(f) - \ve$. By continuity of $f$, compactness of $K$ and local compactness of $X$, there is $O \in \cO(X)$ such that $K \subset O \subset L$ and $f|_O > 1-\ve$. This means that any function $g$ with $g \leq \1_O$ satisfies $\frac f {1-\ve} > g$, and so $\frac 1 {1-\ve} \zeta(f) \geq \tau(O)$. Putting this together, we obtain
$$\tau(K) \geq \zeta(f) - \ve \geq (1-\ve)\tau(O)-\ve \geq \tau(O) - \ve(1+\tau(L)) \geq I-\ve(1+\tau(L))\,.$$
Since $\ve$ was arbitrary and $L$ is fixed, we get $\tau(K) \geq I$, as desired. A similar argument shows inner regularity.

\tb{Additivity}: again, for pairs of disjoint compacts and for pairs of disjoint open sets this is more or less clear, that is follows fairly directly from the definitions and the properties of $\inf$ and $\sup$. The difficulty is to establish additivity for a pair $K\in\cK(X)$, $O \in\cO(X)$, which are disjoint, and such that $K\cup O$ is either open or compact. Let us assume that $U = K\cup O$ is open (and then necessarily with compact closure); the case when the union is compact is treated similarly. Note that regularity implies $\tau(U) \geq \tau(O)+\tau(K)$, since for any compact $K' \subset O$ the union $K'\cup K$ is disjoint, compact and contained in $U$, so $\tau(U) \geq \tau(K) + \tau(K')$. Taking the supremum over all such $K'$, we obtain the statement. Thus it remains to show $\tau(U) \leq \tau(O)+\tau(K)$.

Also by the regularity of $\tau$, we have the following statement: for any $\ve > 0$ there is an open neighborhood $P$ of $K$ with compact closure such that whenever $f$ satisfies $\1_K \leq f \leq \1_P$, it is true that $\zeta(f)\geq\tau(K)\geq\zeta(f)-\ve$. We can choose this $P$ to lie inside any open neighborhood of $K$. Similarly, for any $\ve >0$ there is a compact set $L \subset O$ such that if $g$ satisfies $\1_L \leq g \leq \1_O$, we have $\zeta(g)\leq \tau(O) \leq \zeta(g)+\ve$, and this $L$ can be chosen to contain any prescribed compact subset of $O$.

Let $\ve >0$. Let $L$ be a compact subset of $O$, as we just described. Similarly, let $P$ be an open neighborhood of $K$ with compact closure. We may assume that $\ol P$ is contained in $U-L$. And finally, let $M$ be a compact subset of $U$, containing $L\cup \ol P$, which has the same property with respect to $U$, that is $\1_M\leq h \leq \1_U$ implies $\zeta(h) \leq \tau(U) \leq \zeta(h)+\ve$. Let $h$ be such a function. Let $h'' \fc X \to [0,\ve]$ be such that $h''|_K =\ve$ and $h'' = 0$ outside $P$, and set $h'=h+h''$. Then $h'|_K = 1+\ve$, $1\leq h' \leq 1+\ve$ on $\ol P$ and $h'=h$ outside $P$. Consider two functions $\phi,\psi \fc [0,1+\ve] \to [0,1]$ such that $\phi(t)=0$ for $t\in[0,1]$, $\phi(1+\ve)=1$ and $\psi(t)=t$ for $t\in[0,1]$ and $\phi(t)+\psi(t)=1$ for $t\in[1,1+\ve]$. Define $f=\phi\circ h'$, $g=\psi \circ h'$. These functions have the following properties: $\1_K \leq f \leq \1_P$, $\1_L \leq g \leq \1_O$, $\1_M \leq f+g \leq \1_U$. It follows that
$$\tau(K)+\tau(O) \geq \zeta(f)+\zeta(g)-\ve=\zeta(f+g)-\ve \geq \tau(U)-2\ve\,,$$
where the equality is due to the quasi-linearity of $\zeta$. Thus we obtained the required inequality. \qed
\end{prf}

\subsubsection{From topological measures to quasi-integrals}

Here we fix a topological measure $\tau$ and construct the corresponding quasi-integral $\zeta_\tau$. Let $O \in \cO(X)$ and let $\widehat\tau = \widehat\tau_O$ be the topological measure induced on $\widehat O$ by the one-point compactification procedure. It gives rise to a quasi-integral $\zeta_O \fc C(\widehat O) \to \R$ via the formula
$$\zeta_O(f)=\widehat\tau(\widehat O)\cdot \min f + \int_{\min f}^{\max f} \widehat \tau(\{f \geq t\})\,dt\,.$$
Since $\zeta_O$ is monotone, quasi-linear, and Lipschitz with constant $N_O=\widehat\tau(\widehat O)=\tau(O)$, the same properties are valid for its restriction to $C_c(O) \subset C(\widehat O)$. Therefore, if $f \in C_c(X)$ has support in some $O \in \cO(X)$, define $\zeta_\tau(f)=\zeta_O(f)$. The only thing that we need to check is that this is a correct definition, that is, if the support of $f$ is contained in $O' \in \cO(X)$, then $\zeta_O(f)=\zeta_{O'}(f)$. Since $O\cap O'$ is also in $\cO(X)$ and still contains the support of $f$, we see that it suffices to consider the case $O \subset O'$. Now, both $\zeta_O$ and $\zeta_{O'}$ are quasi-integrals, hence $\zeta_O(f)=\zeta_O(f^+)-\zeta_O(f^-)$, and similarly for $\zeta_{O'}$, where $f^+(x)=\min(0,f(x))$ and $f^-(x)=-\max(0,f(x))$. It follows that we may assume $f \geq 0$. Since $f$ has compact support, $\min f = 0$ on both $\widehat O$ and $\widehat{O'}$; it also has the same maximum. Now for $t > 0$ the set $\{f \geq t\}$ is compact and contained in the support of $f$. By the definition of $\widehat \tau_O$ and $\widehat\tau_{O'}$ we know
$\widehat \tau_O(\{f\geq t\})=\tau(\{f\geq t\})=\widehat\tau_{O'}(\{f\geq t\})$, which means that the two functions coincide on $(0,\max f]$, and hence so do their integrals, which are equal, respectively, to $\zeta_O(f)$, $\zeta_{O'}(f)$.

\subsubsection{The bijection}

Now we have to show that the above procedures of going from quasi-integrals to topological measures and back are inverse to each other.

\begin{prop}\begin{enumerate}
\item Let $\tau$ be a topological measure. Then $\tau_{\zeta_\tau} = \tau$;
\item Let $\zeta$ be a quasi-integral. Then and $\zeta_{\tau_\zeta} = \zeta$.
\end{enumerate}
\end{prop}

\begin{prf}

(i) Let $\sigma = \tau_{\zeta_\tau}$. Recall that if $O \in \cO(X)$, there is a topological measure on $\widehat O$ induced by $\tau$, $\widehat \tau_O$, and that $\zeta_\tau$ restricted to $C_c(O)$ coincides with the restriction of the quasi-integral $\zeta_O$ corresponding to $\widehat \tau_O$. Now, if $K\subset O$ is compact, then
$$\tau(K) = \widehat\tau_O(K)=\inf\{\zeta_O(f)\,|\,f\in C(\widehat O),\,f\geq\1_K\}\,.$$
One can show that the value of the infimum remains unchanged if we only consider functions with compact support in $O$, and then it also equals $\sigma(K)$, by the definition of $\sigma$. Thus $\sigma=\tau$ on $\cK(X)$, and by inner regularity, the same is true on $\cO(X)$.

(ii) Let $\eta = \zeta_{\tau_\zeta}$. Since both $\zeta$ and $\eta$ are quasi-integrals, they respect the positive-negative decomposition of functions, namely $\zeta(f)=\zeta(f^+)-\zeta(f^-)$ and same for $\eta$. Thus it suffices to show that $\zeta$ and $\eta$ coincide on nonnegative functions.

Fix $O \in \cO(X)$ and consider $\widehat\zeta_O$, the quasi-integral induced on $\widehat O$ from $\zeta$ by the one-point compactification procedure. It is represented by a topological measure $\tau'$ on $\widehat O$. Since the restrictions of $\zeta$ and $\widehat \zeta_O$ to $C_c(O)$ coincide, we have $\tau_\zeta(K)=\tau'(K)$ for any compact $K \subset O$. Let $f \in C_c(O)$ be nonnegative. Then
$$\zeta(f)=\widehat\zeta_O(f)=\int_0^{\max f}\tau'(\{f\geq t\})\,dt=\int_0^{\max f}\tau_\zeta(\{f\geq t\})\,dt = \eta(f)\,,$$
the last equality being valid by the definition of $\eta$. \qed
\end{prf}

The proof of Theorem \ref{thm_representation} is thereby complete. \qed

\begin{rem}\label{rem_other_generalizations_qis}Aarnes's representation theorem was generalized to various settings; we refer the reader to \cite{Boardman_qms_cpl_reg_spaces}, \cite{Grubb_LaBerge_spaces_of_qm}, \cite{Wheeler_qms_dim_thry}. The Borel quasi-measures, due to Boardman \cite{Boardman_qms_cpl_reg_spaces}, on a completely regular space are assumed to be defined on all closed and open subsets, and the corresponding quasi-integral is then defined on the space of all bounded continuous functions. As far as we know, quasi-integrals on the space of continuous functions with compact support have not been treated. If the space is assumed to be locally compact, these quasi-integrals are more general than Borel quasi-integrals, in that every Borel quasi-integral determines one as we defined in this paper. Topological measures as defined here are a generalization of both Aarnes topological measures on compact spaces and Radon measures on locally compact spaces. Let us also note that a topological measure $\tau$ on a locally compact space $X$ extends to a unique topological measure $\widehat \tau$ on the one-point compactification $\widehat X$, such that $\widehat \tau(\infty)=0$ if and only if $\tau$ is bounded. This is the case if and only if the corresponding quasi-integral is globally Lipschitz, and the Lipschitz constant evidently equals $\widehat \tau(\widehat X) = \sup_{\cA(X)}\tau$.
\end{rem}

\subsection{Proof of Lemma \ref{thm_Viterbo_fcnls_are_qis}}\label{section_prf_Viterbo_fcnls_are_qis}

Monotonicity follows from that of $\cH$ and of integration. The fact that $\eta_\mu$ is linear on Poisson commuting subspaces of $C^\infty_c(X)$ follows from the `Strong quasi-linearity' property of $\cH$.

It suffices to establish Lipschitz continuity for compact subsets of $T^*\T^n$ of the form $\T^n \times K$ where $K \subset \R^n$ is compact. If a function $f$ has support in $\T^n \times K$, the `Lagrangian' property of $\cH$ implies that $\cH(f)(p) = 0$ for $p\notin K$. Let now $g$ be another function with support in $\T^n \times K$. Then we have
\begin{multline*}
|\eta_\mu(f)-\eta_\mu(g)| \leq \int_{\R^n}|\cH(f)(p)-\cH(g)(p)|\,d\mu(p) \leq \\ \leq \mu(K) \|\cH(f)-\cH(g)\| \leq\mu(K)\|f-g\|\,,
\end{multline*}
which proves Lipschitz continuity together with the announced bound on the Lipschitz constant.

Since integration against a measure is linear, in order to prove the quasi-linearity of $\eta_\mu$, we need only show that $\cH$ is linear on any subspace of $C_c(T^*\T^n)$ of the form $\{\phi\circ f\,|\, \phi \in C(\R),\,\phi(0)=0\}$ for $f\in C_c(T^*\T^n)$. Let $f \in C_c(T^*\T^n)$ and $\phi, \psi \in C(\R)$ with $\phi(0)=\psi(0)=0$; we can replace $\phi$ and $\psi$ by functions with compact support without altering $\phi \circ f$ and $\psi \circ f$. There are functions $f_k \in C^\infty_c(T^*\T^n)$, $\phi_k,\psi_k \in C^\infty_c(\R)$ with $\phi_k(0)=\psi_k(0)=0$ for all $k \in \N$, such that $f_k \to f$, $\phi_k \to \phi$, $\psi_k \to \psi$, all in the $C^0$ norm. It follows that
$$\cH(\phi_k\circ f_k + \psi_k\circ f_k) \to \cH(\phi\circ f + \psi \circ f), \cH(\phi_k\circ f_k) \to \cH(\phi\circ f), \cH(\psi_k\circ f_k) \to \cH(\psi \circ f),$$
all in the $C^0$ norm, due to the Lipschitz continuity of $\cH$. We have
$$\{\phi_k\circ f_k,\psi_k\circ f_k\}=\phi'_k\circ f_k\cdot \psi'_k\circ f_k\{f_k,f_k\}=0$$
for all $k$, which implies, together with the `Strong quasi-linearity' of $\cH$:
$$\cH(\phi_k\circ f_k + \psi_k\circ f_k)= \cH(\phi_k\circ f_k)+ \cH(\psi_k\circ f_k)\,.$$
As $k \to \infty$, the left-hand side of this equality tends to $\cH(\phi\circ f + \psi \circ f)$, while the right-hand side tends to $\cH(\phi\circ f) +\cH(\psi \circ f)$, proving what we wanted. \qed

\section{Proof of the main results and computations}

We use the following result as our main computational tool.
\begin{lemma}\label{lemma_tms_uniquely_detd_by_sbmfds}Let $\tau$ be a topological measure on a manifold without boundary $M$. Then $\tau$ is uniquely determined by its values on codimension $0$ compact connected submanifolds with boundary of $M$.
\end{lemma}
\begin{prf}This result for the case when $M$ is closed was established in \cite{Zapolsky_tms_higher_genus}. Repeated verbatim, that proof also shows that if $M$ is without boundary, the values of $\tau$ on compact subsets of $M$ are uniquely determined by its values on submanifolds as mentioned in the lemma, but without the connectedness assumption. Since a compact manifold with boundary has only finitely many connected components and all of them are also compact subsets of $M$, the additivity of $\tau$ suffices in order to restrict attention to connected submanifolds. Inner regularity allows us to conclude that the values of $\tau$ on $\cO(M)$ are then also uniquely determined. \qed
\end{prf}

\subsection{The topological measure representing $\eta_{p_0}$ for $n=1$}\label{section_computation}
Lemma \ref{thm_Viterbo_fcnls_are_qis} states that if $p_0 \in \R$, then
$$ \eta_{p_0} := \eta_{\delta_{p_0}} = \cH(\cdot)(p_0)\fc C_c(T^*S^1) \to \R$$
is a quasi-integral. According to the representation theorem, Theorem \ref{thm_representation} $\eta_{p_0}$ determines and is determined by a unique topological measure
$$\tau_{p_0} := \tau_{\eta_{p_0}} \fc \cA(T^*S^1) \to \lbrack 0, \infty).$$

Lemma \ref{lemma_tms_uniquely_detd_by_sbmfds} implies that $\tau_{p_0}$ is uniquely determined by its values on compact connected subsurfaces with boundary of $T^*S^1$, and it is these values that we are going to compute.

Such subsurfaces in $T^*S^1$ come in two types. The first type consists of subsurfaces with only contractible boundary components; any such subsurface is a closed disk with holes. Subsurfaces of the second type are those with exactly two non-contractible boundary components; such a subsurface is a closed annulus with holes, the annulus being Hamiltonianly isotopic to a standard one, that is, an annulus of the form $S^1\times[a,b]$.
\begin{lemma}\label{lemma_computation_tau_p} (i) $\tau_{p_0}$ vanishes on subsurfaces of the first type; (ii) if $W$ is of the second type and the annulus is Hamiltonianly isotopic to $S^1\times[a,b]$, then
\begin{equation*}
\tau_{p_0}(W) = {\textup\1}_{[a,b]}(p_0) = \begin{cases} 0 & \text{ if } p_0 \notin \lbrack a,b \rbrack \\
1 &  \text{ if } p_0 \in \lbrack a,b \rbrack
\end{cases}\,.
\end{equation*}
\end{lemma}

\begin{prf}
To prove (i), we note that the homogenization operator $\cH$ vanishes on functions with compact support in a disk.
Indeed, let $U \subset T^*S^1$ be an open disk, let $f$ have support in $U$, and let $p \in \R$. Then there is a Lagrangian $L \subset T^*S^1$ which is Hamiltonianly isotopic to a standard one $S^1 \times \{p\}$, and which avoids $U$. It follows that $f=0$ on $L$. The properties `Invariance' and `Lagrangian' of $\cH$ yield $\cH(f)(p)=0$, and therefore $\cH(f)=0$. From the definition of $\tau_{p_0}$ it follows that $\tau_{p_0}(U)=0$. If $W$ is a subsurface contained in a closed disk, let $V\supset W$ be a slightly larger open disk. Then, using the monotonicity of $\tau_{p_0}$, we obtain
$$0\leq \tau_{p_0}(W) \leq \tau_{p_0}(V)=\sup\{\cH(f)(p_0)\,|\,f\in C_c(T^*S^1),\,f\leq\1_V\}=0\,.$$

Let us turn to (ii). Let $W \subset T^*S^1$ be a subsurface of the second type, that is $W=A-\bigcup_i U_i$, where $A$ is a closed annulus and $U_i \subset A$ are open disks with disjoint closures. Due to the additivity of $\tau_{p_0}$ and the fact that $\tau_{p_0}$ vanishes on open disks, we obtain
$$\tau_{p_0}(W)=\tau_{p_0}(A)-\sum_i\tau_{p_0}(U_i)=\tau_{p_0}(A)\,.$$

By assumption, $A$ is Hamiltonianly isotopic to a standard annulus $S^1 \times \lbrack a,b \rbrack$ and since $\cH$ is invariant under Hamiltonian isotopies it suffices to compute $\tau_{p_0}$ on annuli of this form.

If $p_0 \notin [a,b]$, let $\phi \fc \R \to [0,1]$ be a continuous function such that $\phi(p)=1$ for $p \in [a,b]$ and $\phi(p)=0$ for $p \notin (a-\ve,b+\ve)$, where $\ve > 0$ is chosen small enough so that $p_0 \notin (a-\ve,b+\ve)$. Define $f \in C_c(T^*S^1)$ by $f(q,p)=\phi(p)$. Again, invoking the `Lagrangian' property of $\cH$, we see that $\cH(f)(p_0)=0$. By definition,
$$\tau_{p_0}(S^1\times[a,b])=\inf\{\cH(g)(p_0)\,|\,g\in C_c(T^*S^1),\,g \geq \1_A\}\,.$$
Since $f$ is one of the functions appearing in the latter infimum, it equals $0$, and consequently so does $\tau_{p_0}(S^1\times[a,b])$.

If, on the other hand, $p_0 \in [a,b]$, then for any function $f \in C_c(T^*S^1)$ which equals $1$ on $S^1\times[a,b]$ we will obtain $\cH(f)(p_0) = 1$, once more using property `Lagrangian'. It follows that
$$\tau_{p_0}(S^1\times[a,b])=\inf\{\cH(f)(p_0)\,|\,f\in C_c(T^*S^1),\,f \geq \1_A\} = 1\,.$$
This establishes claim (ii). \qed
\end{prf}

\begin{prf}[of Lemma\ref{lemma_unique_homogenization}]
The existence of the homogenization operator was established by Viterbo \cite{Viterbo_homogenization}. As for uniqueness, it is determined by the functionals $\cH(\cdot)(p)$ for $p\in \R$. As we mentioned above, any such functional is a quasi-integral and thus is uniquely determined by the topological measure $\tau_p$ representing it. In turn, this topological measure is reconstructible from its values on compact connected subsurfaces with boundary of $T^*S^1$. We computed these values above, using only the properties of the homogenization operator stated in the formulation, which implies that it is uniquely determined by those properties. \qed
\end{prf}

\subsection{Proof of Theorem \ref{thm_two_dim}}

Maintain the notations introduced in subsection \ref{section_dim_two}. As we said before, the idea is to compare the topological measures $\tau_0$ and $\tau_r$, representing $\eta_0$ and $\zeta_r$, respectively, namely we will prove that the restriction of $\tau_0$ to $\cA(U_r)$ coincides with $\tau_r$ if and only if $r \in (0,\frac{1}{4} \rbrack$.

The reader should now consult subsection \ref{section_Calabi_qs} for the relevant results about $\tau$, the topological measure representing $\zeta$.

We assume that $\tau_0|_{U_r} = \tau_r$ and show that $r \in (0,\frac 1 4]$. Suppose $r > \frac 1 4$. Let $D \subset U_r$ be a closed disk of area $\geq \frac{1}{2}$. Lemma \ref{lemma_computation_tau_p} states that the topological measures $\tau_p$ vanish on disks, therefore $\tau_0(D)=0$. On the other hand $j(D) \subset S^2$ is a closed disk of area $\geq\frac 1 2$ and thus $\tau_r(D)=\tau(j(D))=1 \neq 0 = \tau_0(D)$, contradiction.

Now we show that $\tau_r=\tau_0|_{U_r}$ for $r \in (0,\frac 1 4]$. Once again, being topological measures on $U_r$, $\tau_0$ and $\tau_r$ are determined by their values on compact subsurfaces with boundary. Therefore it suffices to show that $\tau_0$ and $\tau_r$ coincide on the family of such subsurfaces.

If $W \subset U_r$ is of the first type (see subsection \ref{section_computation}), then $\tau_0(W) = 0$. Now since $W$ is contained in a closed disk $\subset U_r$, the assumption on $r$ implies that this closed disk has area $<\frac 1 2$, and consequently $j(W)$ is contained in a closed disk of area $<\frac 1 2$ as well, implying $\tau_r(W)=\tau(j(W))=0 = \tau_0(W)$.

If $W \subset U_r$ is of the second type, it is an annulus with holes, the annulus being Hamiltonianly isotopic to a standard one $A = S^1\times[a,b]$. By Lemma \ref{lemma_computation_tau_p}, $\tau_0(W)=\tau_0(A)=\1_{[a,b]}(0)$. Let us now compute $\tau_r(W) = \tau(j(W))$. The subsurface $j(W)$ equals $S^2-\big(D_-\cup D_+ \cup \bigcup_iD_i\big)$, where $D_+$ and $D_-$ are the open disks in the complement of $j(W)$, containing the north and the south poles, respectively (we assume without loss of generality that $j(U_r)$ avoids the poles), and $D_i$ are a finite number of open disks with pairwise disjoint closures contained in $j(U_r)$. Since the area of $j(U_r)$ is at most $\frac 1 2$, each $D_i$ has area $\leq \frac 1 2$ and thus $\tau(D_i)=0$, which yields $\tau(j(W))=\tau\big(S^2-(D_-\cup D_+)\big)$. By area arguments it can be seen that $0 \in [a,b]$ if and only if the areas of $D_{\pm}$ are $\leq \frac 1 2$, that is if and only if $\tau\big(S^2-(D_-\cup D_+)\big) = 1$, while $0 \notin [a,b]$ if and only if one of the disks $D_{\pm}$ has area $>\frac 1 2$ if and only if $\tau\big(S^2-(D_-\cup D_+)\big) = 0$. We thus established that $\tau_0(W)=\1_{[a,b]}(0)=\tau(j(W))=\tau_r(W)$, as required. \qed

\subsection{Higher dimensions}\label{section_higher_dims}

In this subsection we prove Proposition \ref{prop_counterex_higher_dim}.

\begin{prf}[of Proposition \ref{prop_counterex_higher_dim}] We denote by $\tau$ the topological measure representing the Calabi quasi-state $\zeta$. The idea of the proof is to construct a Lagrangian torus $L \subset \C P^n$ which is Hamiltonian isotopic to the Clifford torus and which lies in the image $j(U_{\ve,\delta})$, such that $j^{-1}(L)$ equals the product $\ell \times (\T^{n-1} \times \{0\}) \subset T^*S^1 \times T^*\T^{n-1} = T^*\T^n$, where $\ell \subset T^*S^1$ is contractible. Suppose that we constructed such a torus. Then $\tau(L) = \tau(\T^n_{\text{Clif}})=1$, while an argument similar to that of the proof of (i) of Lemma \ref{lemma_computation_tau_p} shows that $\tau_0(j^{-1})(L)=0 \neq 1 = \tau(L)$, showing that $j^*\tau \neq \tau_0$ on $U_{\ve,\delta}$ which implies that $j^*\zeta \neq \eta_0$ on $U_{\ve,\delta}$, therefore proving the desired claim.

We now turn to details. First of all, we give a formula for the symplectic embedding $j \fc U \to \C P^n$ which commutes with the moment maps, as announced in subsection \ref{section_intro_higher_dims}.

For $k \in \N$ we denote by $B^k(r)$ the open Euclidean ball of radius $r$ in $\C^k$, centered at the origin. We normalize the Fubini-Study form $\omega_{\text{FS}}$ on $\C P^n$ by $\int_{\C P^1}\omega_{\text{FS}} = 1$. Consider the symplectic embedding
$$\iota \fc B^n(1) \to \C P^n\,,\iota(z)=\textstyle\Big[\sqrt{1-\sum_j|z_j|^2}:z_1:\dots:z_n\Big]\,,$$
where the symplectic form on $\C^n$ is $\frac 1 \pi$ times the standard one. Recall that
$$\textstyle U:=U_{\delta,\ve}=S^1 \times (-\frac 1{n+1}+\ve,2\ve) \times \big(S^1 \times (-\delta,\delta)\big)^{n-1}\,,$$
and define
$$\kappa \fc U \to B^n(1) \quad\text{by}\quad\textstyle\kappa(q,p) = \big((p_1+\frac 1 {n+1})^{1/2}e^{2\pi i q_1}, \dots, (p_n+\frac 1 {n+1})^{1/2}e^{2\pi i q_n}\big)\,.$$
Set $j = \iota\circ \kappa$. A direct computation shows that $\Psi = \Phi \circ j$. The Clifford torus $\T^n_{\text{Clif}}$ satisfies
$$\iota^{-1}(\T^n_{\text{Clif}}) = \textstyle\{z\in B^n(1)\,|\,|z_j|^2 = \frac 1 {n+1}\,\forall j\}$$
and therefore $j(\T^n \times \{0\}) = \T^n_{\text{Clif}}$.

\begin{figure}[h]
\begin{minipage}{0.48\linewidth}
\centering
\epsfig{file=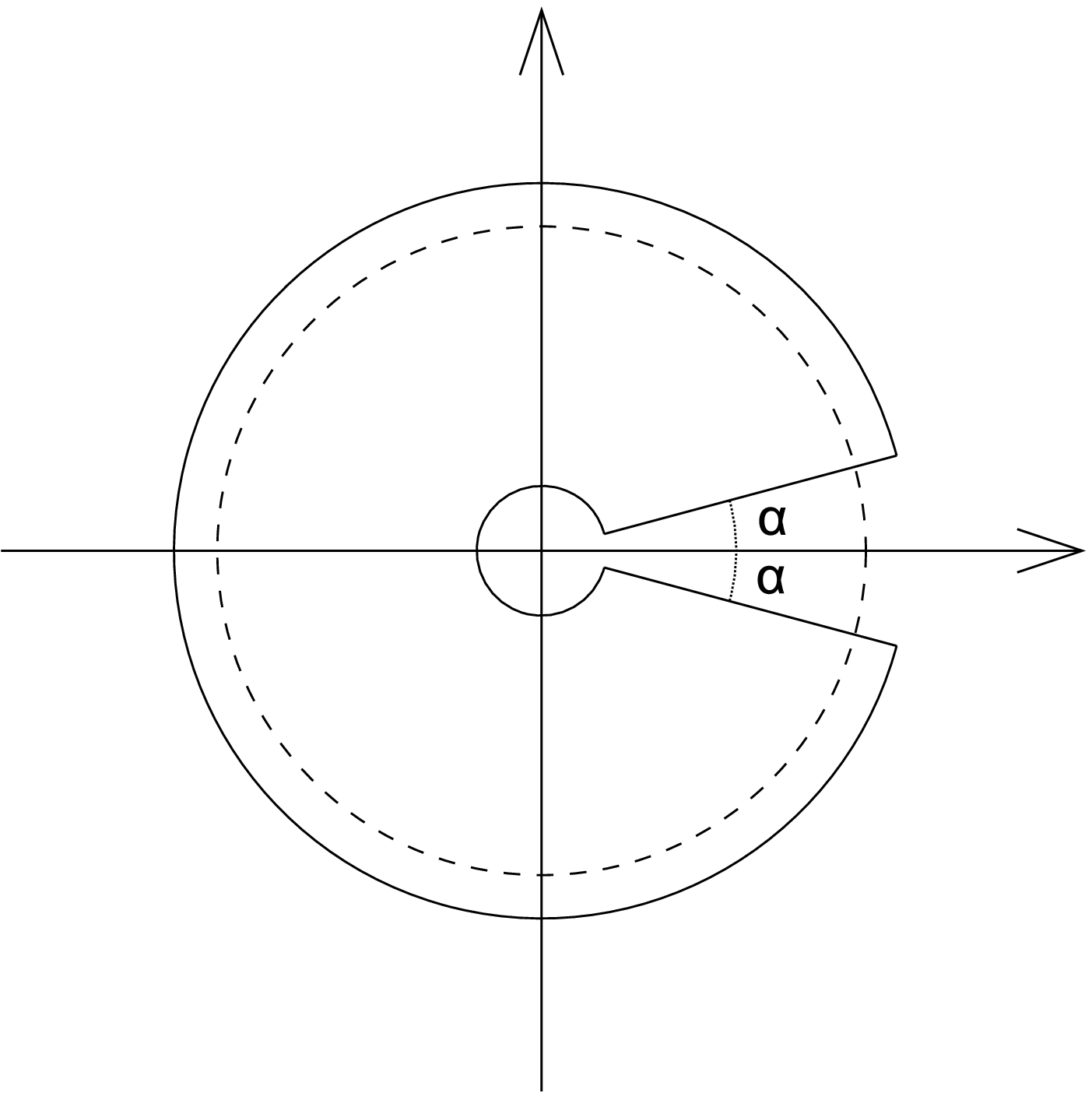,width=6cm}
\caption{The dotted line is the circle $\{|z|^2 = \frac 1 {n+1}\}$; the solid line is the curve $\gamma_0$. The inner and outer radii are $\rho_1,\rho_2$, respectively.}
\label{fig_gamma_zero}
\end{minipage}
\quad
\begin{minipage}{0.45\linewidth}
\centering
\epsfig{file=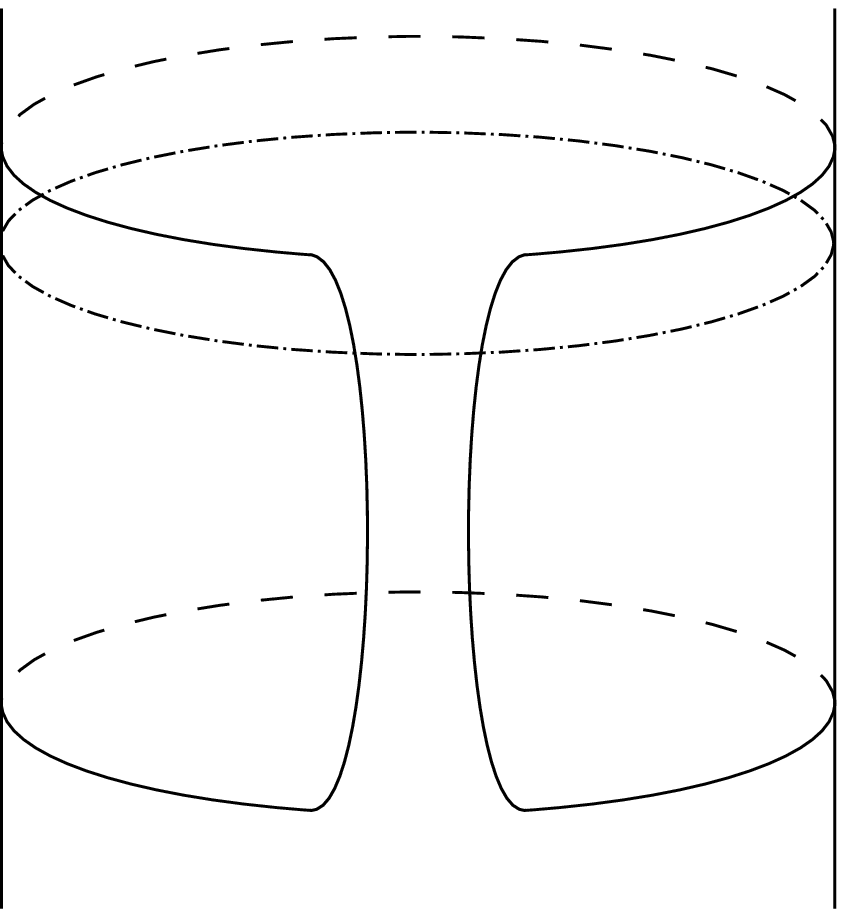,width=5.5cm}
\caption{The dash-dotted curve is the zero section; the solid and dashed curve is $\ell$.}
\label{fig_cylinder}
\end{minipage}
\end{figure}

Let $\rho_1,\rho_2,\alpha > 0$ be real numbers subject to the following conditions: $\frac 1 {n+1} < \rho_2^2 < \frac 1 {n+1} + \ve$, $\alpha < \frac \pi 2$, and $\textstyle(1-\frac \alpha \pi)(\rho_2^2-\frac 1 {n+1}) = \rho_1^2+\frac\alpha\pi(\frac 1 {n+1} - \rho_1^2)$. It is easy to see that such numbers always exist. The last condition expresses the equality of certain areas bounded by two curves, see below. Now consider the following points in the complex plane:
$$\eta_1 = \rho_2e^{ i \alpha},\,\eta_2 = \rho_2e^{- i \alpha},\,\eta_3 = \rho_1e^{- i \alpha},\,\eta_4 = \rho_1e^{ i \alpha}\,,$$
and let $\gamma_0$ be a continuous curve starting at $\eta_1$, continuing counterclockwise along the circle $|z|=\rho_2$ until $\eta_2$, then going to $\eta_3$ along a straight line, following the circle $|z|=\rho_1$ clockwise until $\eta_4$ and finally connecting back to $\eta_1$, again along a straight line, see figure \ref{fig_gamma_zero}. Let $\psi_1$ be a Hamiltonian diffeomorphism of $\C$ generated by a time-dependent Hamiltonian with support in the disk $\{|z|^2 < \frac 1 n\}$, and such that it maps the circle $\{|z|^2=\frac 1 {n+1}\}$ to a smooth curve $\gamma$ which is very close to $\gamma_0$, and such that $\gamma \subset \{0 < |z|^2 < \frac 1 {n+1} + \ve\}$. That such a curve and such a diffeomorphism exist follows from the fact that the signed area between the curve $\gamma_0$ and the circle $\{|z|^2=\frac 1 {n+1}\}$ is zero, as follows from the conditions satisfied by the numbers $\rho_{1,2},\alpha$. Let $\psi$ be the Hamiltonian diffeomorphism of $\big(B^1(\frac 1 {\sqrt n})\big)^n$ given by $\psi = \psi_1 \times \id \times \dots \times \id$. This $\psi$ is also generated by a time-dependent Hamiltonian with compact support in the latter set; the set embeds into $B^n(1)$ and consequently into $\C P^n$ via $\iota$. It follows that $\psi$ can be extended by identity to a Hamiltonian diffeomorphism of $\C P^n$, and let us denote this extended diffeomorphism again by $\psi$. Set $L = \psi(\T^n_{\text{Clif}})$. It can be seen that $L = \iota\big(\gamma \times (\{|z|^2 = \frac 1 {n+1}\})^{n-1}\big)$ and that this torus is contained in $j(U)$ with $j^{-1}(L) = \ell \times (\T^{n-1} \times \{0\})$, where $\ell \subset T^*S^1$ is a contractible curve, see figure \ref{fig_cylinder}. We thus constructed a torus with the required properties. The proof is complete. \qed
\end{prf}

\subsection{The quasi-integral $\eta_{p_0}$ and proof of Lemma \ref{lemma_computation_of_c_pm}}\label{section_eta_p_prf_lemma}

Here we present an explicit formula for $\eta_{p_0}(f)$, Lemma \ref{lemma_computation_of_eta_p}, where $f$ is a sufficiently nicely behaved function, in terms of the Reeb graph of $f$.

Call a function $f \in C^\infty_c(T^*S^1)$ nice if there are numbers $p' < p''$ and $\delta>0$ such that (i) $f(q,p)=0$ for $p \notin(p',p'')$, (ii) $f(q,p)$ is independent of $q$ for $p \in (p',p'+\delta] \cup [p''-\delta,p'')$, (iii) $f$ is generic Morse on $S^1\times (p',p'')$. Note that necessarily $f\neq 0$.

\begin{figure}[h]
%\begin{minipage}{0.48\linewidth}
\centering
\epsfig{file=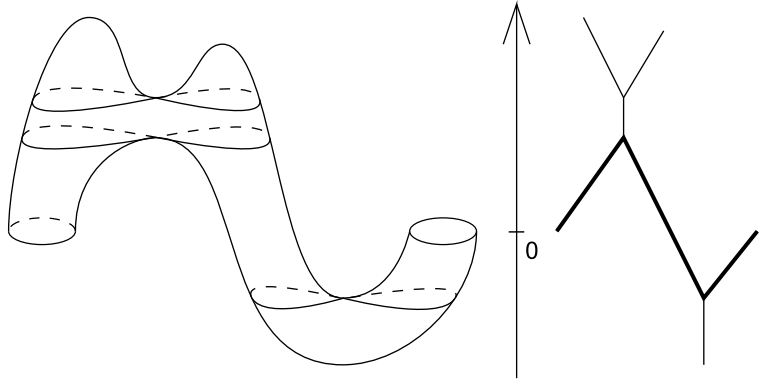,width=10cm}
\caption{The Reeb graph $\Gamma$ (on the right) of the height function on a deformed $T^*S^1$. The bold portion is the subgraph $\Gamma_0$.}
\label{fig_Reeb_graph}
%\end{minipage}
\end{figure}

It is easy to prove
\begin{lemma}The set of nice functions is dense in $C_c(T^*S^1)$ with respect to the $C^0$ norm. \qed
\end{lemma}
\noindent Thus it suffices to compute the values of $\eta_{p_0}$ on nice functions in order to determine it completely. Let therefore $f$ be nice. Define an equivalence relation on $T^*S^1$ by declaring the equivalence classes to be equal to connected components of level sets of $f$. The resulting quotient space is the Reeb graph $\Gamma$ of $f$. It is a tree. Let $\pi \fc T^*S^1 \to \Gamma$ denote the quotient map. Since $f \neq 0$, $\Gamma$ is not a single point. There are two distinguished vertices in $\Gamma$, which we call $v_-$ and $v_+$, which correspond to the connected components $S^1\times (-\infty,p']$ and $S^1\times[p'',\infty)$ of the level set $\{f=0\}$, respectively. Let $\Gamma_0$ be the unique connected linear (that is, of vertex degree $\leq 2$) subgraph containing $v_\pm$. Points in the interior of an edge of $\Gamma_0$ correspond to non-contractible components of regular level sets of $f$, while its vertices other than $v_\pm$ correspond to singular components, that is figures-eight, which are comprised of two non-contractible loops meeting at one point. The reader is referred to figure \ref{fig_Reeb_graph} for an illustration.

We are going to label points of $\Gamma_0$ by subsets of $\R$. To this end, for a non-contractible smooth embedded circle $C \subset T^*S^1$ let $l(C) \in \R$ be its \emph{level}, that is the unique number such that $C$ is Hamiltonianly isotopic to\footnote{The level can be computed as follows. Let $\gamma \fc S^1 \to T^*S^1$ be a parametrization of $C$ such that the composition $\pr\circ\gamma \fc S^1\to S^1$ has degree $1$, where $\pr \fc T^*S^1 \to S^1$ is the projection to the base.  Then $l(C)=\int_{S^1}\gamma^*\lambda$, where $\lambda = p\,dq$ is the Liouville form.} $S^1\times \{l(C)\}$. If $w \in \Gamma_0$ lies in the interior of an edge, label it by the one-point set $\{l(w)\}$, where $l(w) = l(\pi^{-1}(w))$. If $v \in \Gamma_0$ is a vertex other than $v_\pm$, then there are exactly two edges $e_v',e_v''$ of $\Gamma_0$ meeting at it. Let $l_v' =  \lim_{w \in e_v'}l(w)$, $l_v'' =  \lim_{w \in e_v''}l(w)$, where $w$ tends to $v$. These two numbers are never equal\footnote{In fact, their difference is the area bounded by the figure-eight $\pi^{-1}(v)$.} and so we can assume $l_v' < l_v''$ and we label $v$ by the set $[l_v',l_v'']$. Finally, label $v_-$ by $(-\infty,p']$ and $v_+$ by $[p'',\infty)$. Then we have
\begin{lemma}\label{lemma_computation_of_eta_p}$\eta_{p_0}(f)$ equals the value of $f$ at the unique point $w\in\Gamma_0$ such that $p_0$ belongs to the label set of $w$.\qed
\end{lemma}
\noindent A little more abstractly, define a continuous map $\iota \fc \R \to \Gamma_0$ by sending $p \in \R$ to the unique point of $\Gamma_0$ such that $p$ belongs to its label set. There is a continuous function $\ol f \fc \Gamma \to \R$ such that $f = \ol f\circ \pi$. Then $\eta_{p_0}(f) = \ol f (\iota(p_0))$; put differently, $\cH(f)=\ol f \circ \iota$. These claims can be proven using, for example, the techniques of \cite{Zapolsky_Reeb_graph_qs_torus}.

\begin{prf}[of Lemma \ref{lemma_computation_of_c_pm}]We prove that $c_+(f) = \max \cH(f)$. The second statement follows on replacing $f$ by $-f$.

A slightly more delicate version of the `Lagrangian' property of $\cH$, which is proved in \cite{Viterbo_homogenization} implies that if $L \subset T^*S^1$ is a non-contractible smoothly embedded circle, then the existence of a number $c$ such that $f|_L \geq c$ implies $\cH(f)(l(L)) \geq c$, where $l(L)$ is the level of $L$ as above (alternatively, one can use the combination of the `Lagrangian', `Monotonicity' and `Invariance' properties). Since $f|_L \geq \min_L f$, we have $\cH(f)(l(L)) \geq \min_L f$, from which it follows that
$$\max \cH(f) \geq \cH(f)(l(L)) \geq \min_L f\,,$$
and taking supremum over $L \in \cL$ we obtain $\max \cH(f) \geq c_+(f)$.

To prove the reverse inequality we first of all note that both $c_+$ and $\max \cH$ are continuous in the $C^0$ norm, and so it suffices to prove that $c_+(f) \geq \max \cH(f)$ for $f$ a nice function. Choose such an $f$. Note that the function $\ol f \fc \Gamma \to \R$ defined above is strictly monotone on the interiors of the edges of $\Gamma_0$, because interior points correspond to regular components of level sets of $f$. It follows that $\max \cH(f) = \max_{\Gamma_0} \ol f$ is the value of $\ol f$ at one of the vertices $v$ of $\Gamma_0$. Let $e$ be an edge of $\Gamma_0$ adjacent to $v$. The continuity of $\ol f$ implies that $\max \cH(f) = \ol f(v) = \lim_{w\in e} \ol f(w)$, where  the limit is over interior points $w\in e$ tending to $v$. For $w$ an interior point of $e$ the component $\pi^{-1}(w)$ is a smoothly embedded non-contractible circle, and $f$ equals $\ol f(w)$ on it. Therefore $c_+(f) \geq \ol f(w)$. Passing to the limit $w \to v$ we obtain $c_+(f) \geq \lim_{w\in e} \ol f(w) = \max \cH(f)$. \qed
\end{prf}

\end{document}